\documentclass[11pt]{amsart}
\usepackage{bbm,amsfonts,latexsym,amsmath,amscd,amssymb,fancyhdr}
\usepackage[all]{xy}

\numberwithin{equation}{section}

\newcommand{\CC}{\mathbf{C}}
\newcommand{\ZZ}{\mathbf{Z}}

\newcommand{\ii}{\sqrt{-1}}

\newcommand {\be}{\begin{equation}}
\newcommand {\ee}{\end{equation}}

\newcommand{\h}{\begin{eqnarray*}}
\newcommand{\e}{\end{eqnarray*}}

\newtheorem{theorem}{Theorem}[section]

\newtheorem{cor}{Corollary}[section]

\newtheorem{remark}{Remark}[section]

\begin{document}
\title{$\eta$-invariant and modular forms}
\author{Fei Han}
\address{F. Han, Department of Mathematics, National University of Singapore, Block S17, 10 Lower Kent Ridge Road,
Singapore 119076 (mathanf@nus.edu.sg)}
\author{Weiping Zhang}
\address{W. Zhang, Chern Institute of Mathematics \& LPMC, Nankai
University, Tianjin 300071, P.R. China. (weiping@nankai.edu.cn)}

\begin{abstract}
We show that the   Atiyah-Patodi-Singer  reduced $\eta$-invariant of the
twisted Dirac operator on a closed  $4m-1$ dimensional spin manifold,
with the twisted bundle being the Witten bundle appearing in the
theory of elliptic genus, is a meromorphic modular form of weight $2m$ up to an integral
$q$-series. We prove this result by combining our construction of certain
modular characteristic forms associated to a generalized Witten
bundle on spin$^c$-manifolds with a deep
topological theorem due to Hopkins.
\end{abstract}

\maketitle

\section{Introduction and statement of results}

Let $X$ be a smooth manifold. Let $T_{\mathbf C}X$ be the
complexification of the tangent bundle $TX$. One defines the
Witten bundle on $X$ (\cite{W}) as follows, \be
\Theta_q(TX)=\bigotimes_{u=1}^\infty S_{q^u}(T_{\mathbf
C}X-\mathbf{C}^{\mathrm{dim}X}) \otimes \bigotimes_{v=1}^\infty
\Lambda_{-q^{v-{1\over 2}}}(T_{\mathbf
C}X-\mathbf{C}^{\mathrm{dim} X}),\ee where $S_t(\cdot)$ (resp.
$\Lambda_t(\cdot)$) denotes the symmetric  (resp.
exterior) power and $q=e^{2\pi \sqrt{-1} \tau }$ with $\tau \in
\mathbf{H}$, the upper half-plane.

Let $g^{TX}$ be a Riemanian metric on $TX$ and $\nabla^{TX}$ the
associated Levi-Civita connection. If we write \be
\Theta_q(TX)=B_0(TX)+B_1(TX)q^{\frac{1}{2}}+B_2(TX)q+\cdots,\ee
then each $B_i(TX)$ carries a Hermitan metric as well as a
Hermitian connection $\nabla^{B_i(TX)}$ canonically    induced
from $g^{TX}$ and $\nabla^{TX}$. In this way, $\nabla^{TX}$ induces
a Hermitian connection $\nabla^{\Theta_q(TX)}$ on the Witten
bundle $\Theta_q(TX)$.

Now assume that $X$ is closed, spin and of dimension $4m$. Let
$S(TX)=S_+(TX)\oplus S_-(TX)$ be the corresponding Hermitian bundle
of   spinors.
For each $i$, let $D_{X, +}^{B_i(TX)}: \Gamma(S_+(TX)\otimes
B_i(TX))\to \Gamma(S_-(TX)\otimes B_i(TX))$ be the corresponding
twisted  Dirac operator. It is an important and well-known fact
(cf. \cite{Za}) that the $q$-series \be \mathrm{Ind}\left(D_{X,
+}^{\Theta_q(TX)}\right)=\sum_{i=0}^{\infty}\mathrm{Ind}\left(D_{X,
+}^{B_i(TX)}\right)q^{\frac{i}{2}}, \ee which by the Atiyah-Singer
index theorem \cite{AS}  equals to the elliptic genus\footnote{We
refer to \cite[Section 2.1]{HZ2} and \cite[Chapter 1]{Z4} for the
notations of the corresponding characteristic forms appearing
below.}
 \be
 \begin{split}
&\int_X\widehat{A}\left(TX,\nabla^{TX}\right){\rm
ch}\left(\Theta_q(TX),\nabla^{\Theta_q(TX)}\right)\\
=&\sum_{i=0}^{\infty}q^{\frac{i}{2}}
\int_X\widehat{A}\left(TX,\nabla^{TX}\right){\rm
ch}\left(B_i(TX),\nabla^{B_i(TX)}\right),
\end{split}
\ee
is an integral
modular form of weight $2m$ over $\Gamma^0(2)$, where
$\Gamma^0(2)$ is the index 2 modular subgroup of
$SL_2(\mathbf{Z})$ defined by $$
\Gamma^0(2)=\left\{\left.\left(\begin{array}{cc}
a&b\\
c&d
\end{array}\right)\in SL_2(\mathbf{Z})\right|b\equiv0 \ (\rm mod  \
2)\right\}.$$

It is natural to look at what would happen if $X$ is a $4m-1$
dimensional closed spin manifold. In this case, let $E$ be a
Hermitian vector bundle over $X$ carrying a Hermitian  connection
$\nabla^E$. Let $D_X^E:\Gamma(S(TX)\otimes E)\rightarrow
\Gamma(S(TX)\otimes E)$ be the associated twisted Dirac operator,
which is   formally self-adjoint.

Following \cite{APS}, for any ${\rm Re}(s)>>0$, set \be
\eta(D_X^E, s)=\sum_{\lambda\in \mathrm{Spec}(D_X^
E)\setminus\{0\}}\frac{\mathrm{Sgn(\lambda)}}{|\lambda|^s}.\ee By
\cite{APS}, one knows that $ \eta(D_X^E, s)$ is a holomorphic
function in $s$ with $\mathrm{Re}(s)>\frac{\mathrm{dim} X}{2}.$
Moreover, it extends to a meromorphic function over $\mathbf{C}$,
which is holomorphic at $s=0$. The $\eta$ invariant of $D_X^E$, in
the sense of Atiyah-Patodi-Singer \cite{APS}, is defined by
$$\eta(D_X^E)=\eta(D_X^E, 0),$$ while the reduced $\eta$ invariant
is defined and denoted by
$$\overline{\eta}(D_X^E)=\frac{\mathrm{dim}(\mathrm{ker}D_X^E)+\eta(D_X^E)}{2}.$$

It is the aim of this paper to study the modularity of the
$q$-series \be \overline{\eta}\left(D_X^{\Theta_q(TX)}\right)=
\sum_{i=0}^{\infty} \overline{\eta}\left(D_X^
{B_i(TX)}\right)q^{\frac{i}{2}},\ee which is a spectral invariant
depending on   $g^{TX}$.

Assume temporarily that $X$ is the boundary of a $4m$ dimensional
spin manifold $Y$. Let $g^{TY}$ be a Riemannian metric on $TY$
which is of product structure near $\partial Y=X$ and restricts to
$g^{TX}$ on $X$. By the Atiyah-Patodi-Singer index theorem
established in \cite{APS}, one has
\begin{align}\label{0.1}
\int_Y\widehat{A}\left(TY,\nabla^{TY}\right){\rm
ch}\left(\Theta_q(TY),\nabla^{\Theta_q(TY)}\right)-\overline{\eta}\left(D_X^{\Theta_q(TX)}\right)\in{\bf
Z}[[q^{\frac{1}{2}}]].
\end{align}
The term of integration over
$Y$ in (1.7) is a modular form of weight $2m$ over $\Gamma^0(2)$ (similar to the modularity mentioned above for the  elliptic genus in
(1.4),  c.f. \cite{Liu}), although it is not necessary to be an integral modular form
anymore. Therefore, from (\ref{0.1}), one sees that if $X$ bounds
a spin manifold, then for any Riemannian metric on $TX$, $\overline{\eta}
(D_X^{\Theta_q(TX)} )$, up to an integral $q$-series, is a modular
form of weight $2m$ over $\Gamma^0(2)$.

Now the natural question is whether this modularity property for the reduced
$\eta$-invariants holds for any $4m-1$ dimensional closed spin
manifold. The main difficulty of this problem lies in the fact
that given a $4m-1$ dimensional closed spin manifold, it may happen
that it does not bound a spin manifold.

Indeed, it is a well-known fact in cobordism theory that there is
a positive integer $k$ such that $k$ disjoint copies of $X$
bound a spin manifold $\widetilde{Y}$. In this case, one has the
following analogue of (\ref{0.1}),
\begin{align}\label{0.2}
\int_{\widetilde{Y}}\widehat{A}\left(T\widetilde{Y},\nabla^{T\widetilde{Y}}\right){\rm
ch}\left(\Theta_q(T\widetilde{Y}),\nabla^{\Theta_q(T\widetilde{Y})}\right)
-k\,\overline{\eta}\left(D_X^{\Theta_q(TX)}\right)\in{\bf
Z}[[q^{\frac{1}{2}}]].
\end{align} From (\ref{0.2}), one sees that
$\overline{\eta} (D_X^{\Theta_q(TX)} )$ is a modular form up to an
element in $\frac{{\bf Z}[[q^{\frac{1}{2}}]]}{k}$. Thus, the
natural classical method gives the conclusion that
$\overline{\eta}\ (D_X^{\Theta_q(TX)} )$ is a modular form up to
an element in ${\bf Q}[[q^{\frac{1}{2}}]]$ instead of ${\bf
Z}[[q^{\frac{1}{2}}]]$.

On the other hand, if $\widetilde{g}$ is another Riemannian metric
on $TX$ with $\widetilde{\nabla}^{TX}$ being its Levi-Civita
connection and $\widetilde{D}_X^{\Theta_q(TX)}$ being the
corresponding twisted Dirac operator, then by the variation
formula for the reduced $\eta$ invariant (cf. \cite{APS} and
\cite{BF}), one has \be
\overline{\eta}\left(D_X^{\Theta_q(TX)}\right)-\overline{\eta}\left(\widetilde{D}_X^{\Theta_q(TX)}\right)
=\int_X CS_{\Phi}(\nabla^{TX} ,\widetilde{\nabla}^{TX}, \tau)\ \ \
\ \ \mathrm{mod}\,{\bf Z}[[q^{\frac{1}{2}}]],\ee where
$CS_{\Phi}(\nabla^{TX} ,\widetilde{\nabla}^{TX}, \tau)$ is the
Chern-Simons transgression form associated to $\Phi(\nabla^{TX},
\tau)=\left\{\widehat{A}\left(TX,\nabla^{TX}\right){\rm
ch}\left(\Theta_q(TX),\nabla^{\Theta_q(TX)}\right)\right\}^{(4m)}.$
It is easy to see that $\int_XCS_{\Phi}(\nabla^{TX}
,\widetilde{\nabla}^{TX}, \tau)$ is  a modular form of weight $2m$
over $\Gamma^0(2)$ (cf. \cite{CH}). Thus the variation of $\overline{\eta} \left(D_X^{\Theta_q(TX)}\right)$ has mod $\mathbf{Z}$ modularity property. It turns out to be an interesting open
problem that whether $\overline{\eta} \left(D_X^{\Theta_q(TX)}\right)$ is
by itself a modular form of weight $2m$ over $\Gamma^0(2)$ up to an element in
${\bf Z}[[q^{1/2}]]$.

The purpose of this short note is to give an answer to
this question.  Our main result can be stated as follows. 

\begin{theorem}\label{t0.1}
Let $X$ be a  $4m-1$ dimensional  closed spin Riemannian manifold.
Then the reduced $\eta$-invariant $\overline{\eta}
\left(D_X^{\Theta_q(TX)}\right)$ of the twisted Dirac operator
$D_X^{\Theta_q(TX)}   $ is a meromorphic modular form of weight $2m$ over $\Gamma^0(2)$, up to an element in ${\bf
Z}[[q^{\frac{1}{2}}]]$.
\end{theorem}

Here meromorphic modular form is a weaker notion than modular form without requiring holomorphicity but only meromorphicty on the upper half plane.

To prove Theorem \ref{t0.1}, instead of using the cobordism result
as above, we make use of a result due to Hopkins (cf.
\cite[Section 8]{K}) which asserts that for any complex vector
bundle $V$ over $X$, there is a nonnegative integer $s$ such that
 $X\times \CC P^1\times \cdots\times
\CC P^1$ ($s$-copies of $\CC P^1$) bounds a spin manifold $Y$ and
$V\boxtimes H^s$ on $X\times \CC P^1\times \cdots\times \CC
P^1$ extends to $Y$, where $H$ denotes the Hopf hyperplane bundle on $\CC P^1$. We then apply the modular characteristic forms,
which is associated to a generalized Witten bundle
  we have constructed in \cite{HZ2}, on the bounding manifold,
as well as the Atiyah-Patodi-Singer index theorem \cite{APS} to
get the modularity of the reduced $\eta$-invariant in question.

It remains a challenge to find a purely analytic proof of Theorem
\ref{t0.1} without using the deep topological results as
the above mentioned Hopkins' theorem.

Theorem \ref{t0.1} immediately implies that the
quantity  in (\ref{0.2}) is a meromorphic modular form up to an element in
$k{\bf Z}[[q^{\frac{1}{2}}]]$, where $k$ is the positive integer such that $k$ disjoint copies of $X$ bounds $\widetilde{Y}$ as explained before (1.8).   Observe that in (1.8) each $q$-coefficient mod $k$ is a mod $k$ index studied
by Freed and Melrose  in \cite{FM}. It is a topological invariant
and the main result in \cite{FM} provides a topological
interpretation of it.  Therefore, as an application of Theorem 1.1, we have

\begin{cor}\label{t0.2} Let $Y$ be an $4m$ dimensional spin ${\bf Z}/k$-manifold in the sense of Sullivan
(cf. \cite{FM}). Then the mod $k$ index associated to the Witten
bundle $\Theta_q(TY)$ can be represented by a meromorphic modular form
of weight $2m$ over $\Gamma^0(2)$.
\end{cor}

On the other hand, in view of \cite[(25)]{BN}, which corresponds
to the case of $k=1$ in (\ref{0.2}) for the category of stable
almost complex manifolds, Theorem \ref{t0.1} might become a
starting point of a kind of tertiary index theory, in the sense of
\cite[Theorem 4.2]{BN}, for spin manifolds. Recently Ulrich Bunke informed us that Theorem 1.1 can be given an alternative proof by using the theory of the universal $\eta$ invariant (\cite{Bu3}, Lemma 3.1) and a spin version of the $f$-invariant has also been constructed in (\cite{Bu3}, Definition 13.2). 

For completeness, we would like to point out what happens in dimension $4m+1$.  Actually, when $X$ is an $8n+5$ dimensional closed spin manifold, since for each $i$, $ \eta\left(D_X^{B_i(TX)}\right)=0$ and $\mathrm{dim}\left(\mathrm{ker}D_X^{B_i(TX)}\right)$ is even  (c.f. \cite{APS}), we have $\overline{\eta}
\left(D_X^{\Theta_q(TX)}\right)=0\ \mathrm{mod}\, {\bf Z}[[q^{\frac{1}{2}}]]$. In dimension $8n+1$, since $ \eta\left(D_X^{B_i(TX)}\right)=0$  for each $i$ (c.f. \cite{APS}), we have
$\overline{\eta}
\left(D_X^{\Theta_q(TX)}\right)=\frac{\mathrm{dim}\left(\mathrm{ker}D_X^{\Theta_q(TX)}\right)}{2}.$ Therefore in view of the Atiyah-Singer mod 2 index theorem, $\overline{\eta}
\left(D_X^{\Theta_q(TX)}\right)$ can be identified with Ochanine's beta invariant $\beta_q(X)$, the modularity of which has been shown in \cite{O2}.

This paper is organized as follows. In Section 2, we briefly recall
our construction (in \cite{HZ2}) of the modular form
associated to a generalized Witten bundle involving  a complex
line bundle. In Section 3, we combine our modular form and the
Hopkins boundary theorem to prove Theorem 1.1. In Section 4 we
propose a possible refinement of Theorem 1.1 in $8n+3$
dimension.

\section{Complex Line Bundles and Modular Forms} In this section, we briefly review our construction (in \cite{HZ2}) of a modular
form, which is associated to a generalized Witten bundle involving
a complex line bundle.

Let $M$ be a $4l$ dimensional Riemannian manifold. Let
$\nabla^{TM}$ be the associated Levi-Civita connection.

Let $\xi$ be a complex line bundle over $M$. Equivalently, one can view $\xi$ as a rank two real oriented vector bundle over $M$.
 Let $\xi$ carry a Euclidean metric and also a Euclidean connection $\nabla^{\xi}$, let  $c=e(\xi,
\nabla^\xi)$ be  the Euler form   associated to $\nabla^\xi$ (cf.
\cite[Section 3.4]{Z4}). Let $\xi_{\bf C}$ be the complexification
of $\xi$.

If $E$ is a   complex  vector bundle over $M$, set
$\widetilde{E}=E-{\dim E}\in  K(M)$.

Following \cite[(2.5)]{HZ2}, set \be
\begin{split} \Theta_q(TM, \xi)=&\bigotimes_{u=1}^\infty
S_{q^u}(\widetilde{T_\CC M}) \otimes \bigotimes_{v=1}^\infty
\Lambda_{-q^{v-{1\over 2}}}(\widetilde{T_\CC M}-2\widetilde{\xi_\CC})\\
&\otimes\bigotimes_{r=1}^{\infty}\Lambda_{q^{r-\frac{1}{2}}}(\widetilde{\xi_\CC})
\otimes\bigotimes_{t=1}^{\infty}\Lambda_{q^{t}}(\widetilde{\xi_\CC}) ,\\
\end{split}\ee which is
an element in $K(M)[[q^{1\over2}]]$. As before, $\nabla^{TM}$ and
$\nabla^\xi$ induce a Hermitian connection $\nabla^{ \Theta_q(TM,
\xi)}$ on $\Theta_q(TM, \xi)$.

Let $P(TM, \xi, \tau)\in \Omega^{4l}(M)$ be the  characteristic
form defined by
 \be P(TM, \xi, \tau):=\left\{\widehat{A}(TM,
\nabla^{TM})\cosh\left(\frac{c}{2}\right)\mathrm{ch}\left(\Theta_q(TM,
\xi), \nabla^{ \Theta_q(TM, \xi)}\right)\right\}^{(4l)}.\ee

It is shown in \cite{HZ2} that $ P(TM, \xi, \tau)$ can be
expressed by using the formal Chern roots of $(T_\CC M,
\nabla^{T_\CC M})$ and $c$ through the Jacobi theta functions,
which are defined as follows (cf. \cite{C} and \cite[Section
2.3]{HZ2}): \h \theta(v,\tau)=2q^{1/8}\sin(\pi v)
\prod_{j=1}^\infty\left[(1-q^j)(1-e^{2\pi
\sqrt{-1}v}q^j)(1-e^{-2\pi \sqrt{-1}v}q^j)\right]\ ,\e \h
\theta_1(v,\tau)=2q^{1/8}\cos(\pi v)
 \prod_{j=1}^\infty\left[(1-q^j)(1+e^{2\pi \sqrt{-1}v}q^j)
 (1+e^{-2\pi \sqrt{-1}v}q^j)\right]\ ,\e
\h \theta_2(v,\tau)=\prod_{j=1}^\infty\left[(1-q^j)
 (1-e^{2\pi \sqrt{-1}v}q^{j-1/2})(1-e^{-2\pi \sqrt{-1}v}q^{j-1/2})\right]\
 ,\e
\h \theta_3(v,\tau)=\prod_{j=1}^\infty\left[(1-q^j) (1+e^{2\pi
\sqrt{-1}v}q^{j-1/2})(1+e^{-2\pi \sqrt{-1}v}q^{j-1/2})\right]\ .\e
The theta functions are all holomorphic functions for $(v,\tau)\in
\mathbf{C \times H}$, where $\mathbf{C}$ is the complex plane and
$\mathbf{H}$ is the upper half plane. Let  $\{\pm 2\pi \ii x_i\}$
be the
 formal Chern roots for
 $(T_\CC M, \nabla^{T_\CC M})$ and $c=2\pi \ii u$, we have
\be P(TM, \xi, \tau)\\
=\left\{\left(\prod_{i=1}^{2l}x_i\frac{\theta'(0,\tau)}{\theta(x_i,\tau)}
\frac{\theta_{2}(x_i,\tau)}{\theta_{2}(0,\tau)}\right) \frac{\theta_1(u,\tau)}{\theta_1(0,\tau)}\frac{\theta_2^2(0,\tau)}{\theta_2^2(u,\tau)}\frac{\theta_3(u,\tau)}{\theta_3(0,\tau)}\right\}^{(4l)}.\ee

By using the transformation laws of theta functions (cf. \cite{C}
and \cite[Section 2.3]{HZ2}), one sees as in \cite[Proposition
2.6]{HZ2} that $P(TM, \xi, \tau)$ is a modular form of weight $2l$
over $\Gamma^0(2)$.

\section{Proof of the Main Theorem}
In this section, we will prove our main result Theorem 1.1.

The topological tool we will use is the following boundary theorem
of Hopkins (cf. \cite[Section 8]{K}).
\begin{theorem} [Hopkins] Let $X$ be a compact, odd dimensional
spin manifold and $V\to X$ a complex vector bundle over $X$. Then there is an
integer $s$ such that the vector bundle $V\boxtimes (\boxtimes_{j=1}^sH_j) \to X\times (\CC P^1)^s$ is a boundary, where $H_j$ denotes the Hopf hyperplane bundle on the $j$-th copy
of $\CC P^1$. In other words, there is a spin
manifold $Y$ with a complex vector bundle $W$ on
$Y$ such that $W|_{\partial Y}=V\boxtimes (\boxtimes_{j=1}^sH_j),$

\end{theorem}

In what follows, we will combine this Hopkins boundary theorem with
the modular characteristic  form constructed in Section 2 to give
a proof of Theorem 1.1.

$ $

\noindent {\it Proof of Theorem 1.1:}  Without loss of generality, for the $4m-1$ dimensional closed spin manifold $X$, in view of
the Hopkins boundary theorem, we take an even integer $s$ so that
the complex line bundle
$$p^*(\boxtimes_{j=1}^sH_j)\to X\times (\CC P^1)^s$$ bounds, where $p: X\times (\CC P^1)^s \to (\CC P^1)^s$ is the natural projection.  This means that there is a spin manifold $Y$ and a complex line bundle $\zeta$ over $Y$ such that $\partial Y=X\times (\CC P^1)^s$ and
$\zeta|_{X\times (\CC P^1)^s}=p^*(\boxtimes_{j=1}^sH_j).$

Let $g^{TX}$ be any Riemmnian metric on $X$. Equip $\CC P^1$'s
with arbitrary Riemannian metrics and the $H_j$'s with arbitrary
Euclidean metrics and Euclidean connections.

Let $g^{TY}$ be a metric  on $TY$ such that it is of product
structure near $X\times (\CC P^1)^s $ and restricts to the product
metric on $X\times (\CC P^1)^s $. Let $\nabla^{TY}$ be the
Levi-Civita connection associated to $g^{TY}$.

Let $g^{\zeta}$ be  an Euclidean metric on $\zeta$ (viewed as an
oriented real plane bundle) such that
 $g^{\zeta}$   is of product structure
  near $X\times (\CC P^1)^s $  and restricts to the Euclidean metric on $p^*(\boxtimes_{j=1}^sH_j)$ on $X\times (\CC P^1)^s $. Let $\nabla^{\zeta}$ be an
  Euclidean connection of $g^{\zeta}$ which is of product structure  near $X\times (\CC P^1)^s $ and
  restricts to the canonically induced Euclidean connection on $p^*(\boxtimes_{j=1}^sH_j)$ on $X\times (\CC P^1)^s $.

Let $c=e(\zeta)$ and $z_j=\frac{c_1(H_j)}{\pi \ii}, 1\leq j\leq s$.

By applying the Atiyah-Patodi-Singer index theorem \cite{APS} to
the twisted Dirac operator $D_Y^{\Theta_q(TY, \zeta^2)\otimes
\zeta}$, in noting that $$\left.\left( {\Theta_q(TY,
\zeta^2)\otimes \zeta}\right)\right|_{X\times (\CC
P^1)^s}=\Theta_q\left(T\left(X\times (\CC P^1)^s\right),
\left(p^*(\boxtimes_{j=1}^sH_j)\right)^2\right)\otimes
p^*\left(\boxtimes_{j=1}^sH_j\right),$$ one finds that there exist
integers $a_i$'s such that \be
\begin{split}
&\overline{\eta}\left(D_{X\times (\CC P^1)^s}^{\Theta_q(T(X\times (\CC P^1)^s), (p^*(\boxtimes_{j=1}^sH_j))^2)
\otimes p^*(\boxtimes_{j=1}^sH_j)}\right)\\
=&\int_Y \widehat{A}\left(TY,\nabla^{TY}\right){\rm
ch}\left(\Theta_q(TY, \zeta^2)\otimes \zeta,\nabla^{\Theta_q(TY, \zeta^2)\otimes \zeta}\right)-\sum_{i=0}^\infty a_iq^{\frac{i}{2}}\\
=&\int_Y \widehat{A}\left(TY,\nabla^{TY}\right)e^c\,{\rm
ch}\left(\Theta_q(TY, \zeta^2),\nabla^{\Theta_q(TY, \zeta^2)}\right)-\sum_{i=0}^\infty a_iq^{\frac{i}{2}}\\
=&\int_Y \widehat{A}\left(TY,\nabla^{TY}\right)\cosh(c)\,{\rm
ch}\left(\Theta_q(TY, \zeta^2),\nabla^{\Theta_q(TY,
\zeta^2)}\right)-\sum_{i=0}^\infty a_iq^{\frac{i}{2}},\end{split}
\ee where the last equality follows from the fact that  $s$ is an
even integer.

Let $r: X\times (\CC P^1)^s\to X$ be the natural  projection. For bundles $E\to X$ and $F\to (\CC P^1)^s$, by
separation of variables,  we have
$$ \eta\left(D_{X\times (\CC P^1)^s}^{(r^*E)\otimes (p^*F)}\right)=\eta(D_{X}^E)\cdot \mathrm{Ind}(D_{(\CC P^1)^s,+}^F).$$
So we have
$$\overline{\eta}\left(D_{X\times (\CC P^1)^s}^{(r^*E)\otimes (p^*F)}\right)=\overline{\eta}(D_{X}^E)\cdot \mathrm{Ind}(D_{(\CC P^1)^s, +}^F)
+\mathrm{dim} (\mathrm{ker}D_{X}^E)\mathrm{dim}(\mathrm{ker}(D_{(\CC P^1)^s, -}^F)).$$
From the above formula,  we  can see that there are  integers $b_i$'s
such that \be
\begin{split} &\overline{\eta}\left(D_{X\times (\CC P^1)^s}^{\Theta_q(T(X\times (\CC P^1)^s), (p^*(\boxtimes_{j=1}^sH_j))^2)
\otimes p^*(\boxtimes_{j=1}^sH_j)}\right)-\sum_{i=0}^\infty b_iq^{\frac{i}{2}}\\
=& \overline{\eta}\left(D_{X\times (\CC P^1)^s}^{\Theta_q(r^*TX\oplus p^*T(\CC P^1)^s,
(p^*(\boxtimes_{j=1}^sH_j))^2)\otimes p^*(\boxtimes_{j=1}^sH_j)}\right)-\sum_{i=0}^\infty b_iq^{\frac{i}{2}}\\
=&\overline{\eta}\left(D_{X\times (\CC P^1)^s}^{r^*\Theta_q(TX)\otimes p^*(\Theta_q(T(\CC P^1)^s,
(\boxtimes_{j=1}^sH_j)^2)\otimes \boxtimes_{j=1}^sH_j)}\right)-\sum_{i=0}^\infty b_iq^{\frac{i}{2}}\\
=& \overline{\eta}\left(D_X^{\Theta_q(TX)}\right)\cdot \mathrm{Ind}\left(D_{(\CC P^1)^s, +}^{\Theta_q(T(\CC P^1)^s, (\boxtimes_{j=1}^sH_j)^2)\otimes \boxtimes_{j=1}^sH_j}\right)\\
=&\overline{\eta}\left(D_X^{\Theta_q(TX)}\right)\\
&\cdot \int_{(\CC P^1)^s} \widehat{A}(T(\CC P^1)^s, \nabla^{T(\CC P^1)^s})
e^{c_1(H_1)+\cdots+c_1(H_s)}\mathrm{ch}\left(\Theta_q(T(\CC P^1)^s, (\boxtimes_{j=1}^sH_j)^2)\right)\\
=&\overline{\eta}\left(D_X^{\Theta_q(TX)}\right)\\
&\cdot \int_{(\CC P^1)^s} \left( \prod_{j=1}^{s}z_j \frac{\theta'(0,\tau)}{\theta(z_j,\tau)}\frac{\theta_2(z_j,\tau)}{\theta_2(0,\tau)}\right)\frac{\theta_1(\sum_{j=1}^{s}z_j,\tau)}{\theta_1(0,\tau)}\frac{\theta_2^2(0,\tau)}{\theta_2^2(\sum_{j=1}^{s}z_j,\tau)}\frac{\theta_3(\sum_{j=1}^{s}z_j,\tau)}{\theta_3(0,\tau)}\\
=&\overline{\eta}\left(D_X^{\Theta_q(TX)}\right)\cdot \int_{(\CC
P^1)^s}
\frac{\theta_1(\sum_{j=1}^{s}z_j,\tau)}{\theta_1(0,\tau)}\frac{\theta_2^2(0,\tau)}{\theta_2^2(\sum_{j=1}^{s}z_j,\tau)}\frac{\theta_3(\sum_{j=1}^{s}z_j,\tau)}{\theta_3(0,\tau)},\end{split}\ee
where the last equality holds due to the fact that
$\frac{x}{\theta(x,\tau)}$ and $\theta_2(x,\tau)$ are both even
functions about $x$ and $\int_{\CC P^1}z_j^n=0$ if $n>1$.

 Since $s$ is an even integer, from the knowledge about the modular form $P(TM, \xi, \tau)$ constructed
  in Section 2,  we know that
  $$f_s(\tau):=\int_{(\CC P^1)^s} \frac{\theta_1(\sum_{j=1}^{s}z_j,\tau)}{\theta_1(0,\tau)}\frac{\theta_2^2(0,\tau)}
  {\theta_2^2(\sum_{j=1}^{s}z_j,\tau)}\frac{\theta_3(\sum_{j=1}^{s}z_j,\tau)}{\theta_3(0,\tau)}$$ is an integral
  modular form of weight $s$ over $\Gamma^0(2)$. Moreover, since
  $$\int_{(\CC P^1)^s} \widehat{A}\left(T(\CC P^1)^s, \nabla^{T(\CC P^1)^s}\right)
e^{c_1(H_1)+\cdots+c_1(H_s)}=1,$$ we see that $f_s(\tau)$ has constant term $1$. Therefore $f_s^{-1}(\tau)\in {\bf
Z}[[q^{\frac{1}{2}}]]$.

From (3.1) and (3.2),  we have
\be
\begin{split}
&\overline{\eta}\left(D_X^{\Theta_q(TX)}\right)\\
=&f_s^{-1}(\tau)\cdot \int_Y
\widehat{A}\left(TY,\nabla^{TY}\right)\cosh(c)\,{\rm
ch}\left(\Theta_q(TY, \zeta^2),\nabla^{\Theta_q(TY,
\zeta^2)}\right)\\
&-f_s^{-1}(\tau)\cdot\left(\sum_{i=0}^\infty
\left(a_i+b_i\right)q^{\frac{i}{2}}\right).
\end{split}
\ee

Still by the modularity of $P(TM, \xi, \tau)$ constructed in
Section 2, we know that
$$\int_Y \widehat{A}\left(TY,\nabla^{TY}\right)\cosh(c)\,{\rm
ch}\left(\Theta_q(TY, \zeta^2),\nabla^{\Theta_q(TY,
\zeta^2)}\right)$$ is a modular form of weight $2m+s$ over
$\Gamma^0(2)$. So
$$f_s^{-1}(\tau)\cdot \int_Y \widehat{A}\left(TY,\nabla^{TY}\right)\cosh(c)\,{\rm
ch}\left(\Theta_q(TY, \zeta^2),\nabla^{\Theta_q(TY, \zeta^2)}\right)$$ is a meromorphic modular form of weight $2m$ over $\Gamma^0(2)$.

Therefore from (3.3), we see that
 $$\overline{\eta}
(D_X^{\Theta_q(TX)})=f_s^{-1}(\tau)\cdot \int_Y \widehat{A}\left(TY,\nabla^{TY}\right)\cosh(c)\,{\rm
ch}\left(\Theta_q(TY, \zeta^2),\nabla^{\Theta_q(TY, \zeta^2)}\right)\ \ \mathrm{mod}\, {\bf Z}[[q^{\frac{1}{2}}]],$$ a meromorphic modular form of weight $2m$ over
$\Gamma^0(2)$.
The proof of Theorem 1.1 is complete. \ \ \  Q.E.D.

\begin{remark} The modular form $f_s(\tau)$ in the above proof can be explicitly expressed by theta functions and their derivatives. For example, we have
\be
f_2(\tau)=-\frac{1}{\pi^2}\left(\frac{\theta_1''(0,\tau)}{\theta_1(0,\tau)}-2\frac{\theta_2''(0,\tau)}{\theta_2(0,\tau)}+\frac{\theta_3''(0,\tau)}{\theta_3(0,\tau)}\right)\ee
and
\be
\begin{split}
f_4(\tau)=\frac{1}{\pi^4}&\left(\frac{\theta_1^{(4)}(0,\tau)}{\theta_1(0,\tau)}-2\frac{\theta_2^{(4)}(0,\tau)}{\theta_2(0,\tau)}+\frac{\theta_3^{(4)}(0,\tau)}{\theta_3(0,\tau)} +18\left(\frac{\theta_2''(0,\tau)}{\theta_2(0,\tau)}\right)^2
\right. \\
&\ \ \ \left.
-12\frac{\theta_1''(0,\tau)}{\theta_1(0,\tau)}\frac{\theta_2''(0,\tau)}{\theta_2(0,\tau)}
-12\frac{\theta_3''(0,\tau)}{\theta_3(0,\tau)}\frac{\theta_2''(0,\tau)}{\theta_2(0,\tau)}
+6\frac{\theta_1''(0,\tau)}{\theta_1(0,\tau)}\frac{\theta_3''(0,\tau)}{\theta_3(0,\tau)}\right).
\end{split}
\ee
\end{remark}

\begin{remark} Let $X$ be a compact, odd dimensional spin manifold. Define
$$H(X):=\{h\in \ZZ: the\  line\  bundle\  p^*(\boxtimes_{j=1}^hH_j)\to X\times (\CC P^1)^h\ bounds\},$$
where $p: X\times (\CC P^1)^h \to (\CC P^1)^h$ is the natural projection and $H_j$ denotes the Hopf hyperplane bundle on the $j$-th copy
of $\CC P^1$. Define the Hopkins' index of $X$, $h(X):=\min H(X)$. Obviously, when $X$ is a boundary by itself, $h(X)=0$. It is clear that $H(X)=\{s\in \ZZ: s\geq h(X)\}.$ 

In the proof of Theorem 1.1, we may take 
any even number $s\in H(X)$ and denote the corresponding $Y$ and $\zeta$ by $Y_s$ and $\zeta_s$. Then the proof of Theorem 1.1 tells us that, up to an element in ${\bf Z}[[q^{\frac{1}{2}}]]$,
$$\overline{\eta}
(D_X^{\Theta_q(TX)})=f_s^{-1}(\tau)\cdot \int_{Y_s} \widehat{A}\left(TY_s,\nabla^{TY_s}\right)\cosh(e(\xi_s))\,{\rm
ch}\left(\Theta_q(TY_s, \zeta_s^2),\nabla^{\Theta_q(TY_s, \zeta_s^2)}\right).$$
Clearly, if $h(X)=0$, one gets (\ref{0.1}).
Therefore, for every even number $s\geq 2[\frac{h(X)+1}{2}]$, one can construct a meromorphic modular form of weight 2m over
$\Gamma^0(2)$ of above form, that is equal to $\overline{\eta}
(D_X^{\Theta_q(TX)})$ up to an element in ${\bf Z}[[q^{\frac{1}{2}}]]$. The poles of these meromorphic modular forms are just the zeros of the modular forms $f_s(\tau)$. We hope that further study of the modular forms $f_s(\tau)$ will bring better understanding of modularity of $\overline{\eta}
(D_X^{\Theta_q(TX)}).$  
\end{remark}
\begin{remark} We refer to \cite{Bu3} for an alternative approach to the modularity of $\overline{\eta}
(D_X^{\Theta_q(TX)})$, which is shown to be not only a meromorphic modular form but also a modular form using the theory of universal $\eta$-invariant. 
\end{remark}

\section{The cases of dimension $8n+3$}

In this section, we discuss the case of dimension $ 8n+3$. In this dimension, it is known that $\overline{\eta}   (D_X^{\Theta_q(TX)}  )$ is mod
$2{\bf Z}[[q^{\frac{1}{2}}]]$ smooth. That is, in the right hand
side of (1.9), the term mod ${\bf Z}[[q^{\frac{1}{2}}]]$ can be
replaced by  mod $2{\bf Z}[[q^{\frac{1}{2}}]]$. Therefore
it is natural to propose the following conjecture whose statement
refines Theorem 1.1 in this case.

$\ $

\noindent{\bf Conjecture  4.1.} {\it Let $X$ be an $8n+3$
dimensional closed spin Riemannian manifold. Then the reduced
$\eta$-invariant  $\overline{\eta}  (D_X^{\Theta_q(TX)} )$ of the
twisted Dirac operator $D_X^{\Theta_q(TX)}   $ is a meromorphic modular form
of weight $4n+{2}$ over $\Gamma^0(2)$, up to an element in $2{\bf
Z}[[q^{\frac{1}{2}}]]$.}

$\ $

Recall that a mod $2k$ refinement of the Freed-Melrose mod $k$
index  for real vector bundles over $8n+4$ dimensional manifolds
has been defined in \cite[Section 3]{Z2}. In view of this, one can
propose a refinement of Corollary \ref{t0.2}, in the case of $\dim
Y=8n+4$, as follows.

$\ $

\noindent{\bf Conjecture  4.2.} {\it Let $Y$ be an $8n+4$
dimensional spin ${\bf Z}/k$-manifold in the sense of Sullivan
(cf. \cite{FM}). Then the mod $2k$ index associated to the Witten
bundle $\Theta_q(TY)$ can be represented by a meromorphic modular form of
weight $4n+{2}$ over $\Gamma^0(2)$.}

$\ $

By the method of this paper, in order to prove Conjectures  4.1
and 4.2,  one perhaps needs a kind of Hopkins boundary theorem for
real vector bundles. Or, one may try to develop a direct analytic
approach, which, even for Theorem 1.1, is a challenging problem as
we indicated in Section 1.

$\ $

$$ $$

\noindent {\bf Acknowledgement} This work was started  when both
of us were participating an international conference  at UC Santa
Barbara organized by Xianzhe Dai in July, 2012. We would like to
thank Xianzhe Dai for the invitation as well as the kind
hospitality. We would also like to thank Siye Wu for bringing our
attention to \cite{K}.  We are grateful to Ulrich Bunke for the discussion on the topic. 

The work of F. H. was partially supported
by a start-up grant and AcRF R-146-000-163-112 from National University of Singapore. The
work of W. Z. was partially supported by NNSFC and MOEC.

\end{document}